\documentclass[12pt]{article}

\usepackage{amsmath, amsthm, amsfonts, graphicx}
\usepackage[all, cmtip]{xy}

\theoremstyle{plain}
\newtheorem{theorem}{Theorem}
\newtheorem{proposition}[theorem]{Proposition}
\newtheorem{corollary}[theorem]{Corollary}

\theoremstyle{definition}
\newtheorem{axiom}{Axiom}

\newcommand{\cC}{\mathcal{C}}

\newcommand{\Rp}{\mathbb{R}_{\geq 0}}
\newcommand{\e}{\hat{\times}}
\newcommand{\Tr}{\mathop{\mathrm{Tr}}\nolimits}
\newcommand{\id}{\mathop{\mathrm{id}}\nolimits}

\begin{document}

\title{The Category-Theoretic Arithmetic of Information}

\author{Benjamin Allen\\Boston University\\Department of Mathematics and Statistics\\111 Cummington St\\Boston, MA 02115}

\maketitle

\begin{abstract}
We highlight the underlying category-theoretic structure of measures of information flow.  We present an axiomatic framework
in which communication systems are represented as morphisms, and information flow is characterized by its behavior when 
communication systems are combined.  Our framework includes a variety of discrete, continuous, and, conjecturally, quantum information 
measures.  It also includes some familiar mathematical constructs not normally associated with information, such as 
vector space dimension.  We discuss these examples and prove basic results from the axioms.
\end{abstract}

\section{Introduction}

Information theory characterizes the transmission of information using a variety of measures, including discrete and continuous
versions of entropy, mutual information, and channel capacity.  Quantum information theory has added to this list by introducing notions
of quantum channel capacity.  

In this work we identify a category-theoretic structure underlying these measures.  We show that communication systems have a natural
representation as morphisms, and using this representation, we give axioms that characterize information flow.  
Unlike other axiomatizations of information \cite{Shannon, Khinchin, Lee, Abe, Suyari}, we do not require specific mathematical data
such as a random variable or set partition.  
Our framework thereby encompasses not only a single information measure or family of measures, but a wide variety of  
functions used in discrete, continuous, and, conjecturally, quantum information theory.  Vector space 
dimension also satisfies the axioms of an information measure, giving additional support to connections between information and 
dimension that have been noticed elsewhere \cite{Lutz, Hitchcock}.  We hope our work will help unify and extend the results of 
information theory by stimulating the discovery of general theorems and new ways to measure information

The organization of this paper is as follows: In section 2 we introduce the representation of communication systems as morphisms,
and in section 3 we discuss how such systems may be formally combined using category-theoretic operations.  In section 4 we present
our axiomatic framework describing the behavior of information under these operations.  Section 5 explores the various settings to 
which our framework applies, and in section 6 we use the axioms to prove basic results about the general properties of information.

\section{Communication Systems}

\begin{figure}
\includegraphics{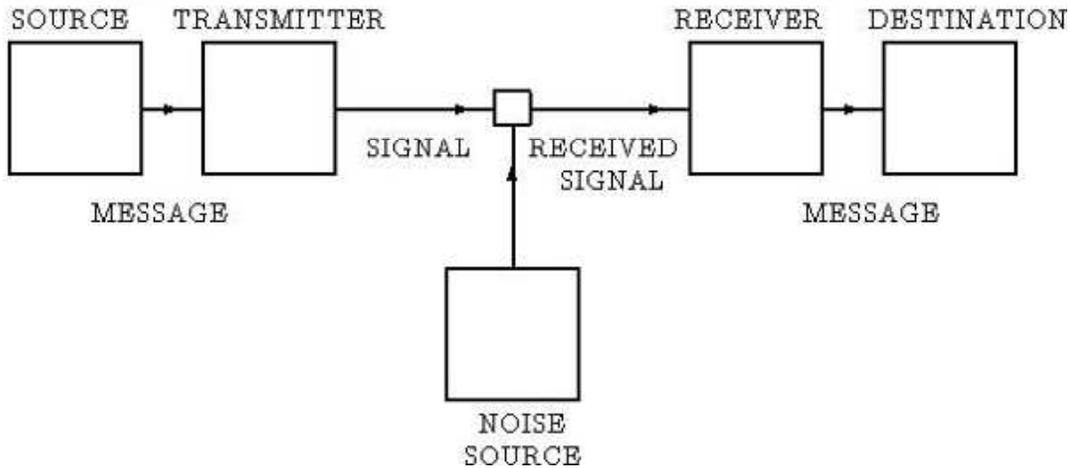}
\caption{Shannon's model of a communcation system \cite{Shannon}.}
\label{ShannonModel}
\end{figure}

We view information as something transmitted by a communication system.  A communication system consists of a source and a 
destination, together with a method of transmitting information between them in the form of messages.  This picture
can be seen as a simplified version of Shannon's model, depicted in figure \ref{ShannonModel}.  

We use the term ``communication system'' very generally.  A person looking at a picture is an example of such a system, since
information is transmitted from the picture to the person.  Our focus here is on systems with mathematical representations, but
the framework we present may also give insight into less rigorous situations.

Communication systems have a natural representation as category-theoretic morphisms.  They are directed relations
between two objects satisfying the three basic properties of morphisms: a) two communication systems can be ``composed'' by passing 
a message through 
one and then the other, b) this composition is associative, and c) ``identity'' communication systems exist, wherein 
the source is the same as the destination and the message stays as is.  We therefore model communication systems as morphisms
in various categories

In most of the categories we consider here, the objects are sets with perhaps some additional structure, and the morphisms
are set mappings preserving the structure.  The elements of the sets represent
possible messages, and the morphisms map messages sent onto messages received.  Information can then be quantified in terms of 
how much the received message tells the receiver about which message was sent.  We note, however, that our formalism may
also be used with more abstract categories (see section \ref{Dual Categories}), wherein the interpretation of morphisms as
communication systems is less intuitive.

\section{Combining Morphisms}

We will characterize information in terms of its behavior when communication systems are combined.  However,
there are several meaningful ways to combine communication systems, each of which corresponds to an operation on morphisms.  
One such operation is composition, as described above.  In this section we define two new operations on morphisms and interpret them
in terms of communication systems.  

Let $\cC$ be a category.  For morphisms $f, g \in \hom(\cC)$, we define the \emph{external product} $f \e g \in \hom{(\cC)}$ to be the 
product, if it exists, of $f$ and $g$ in the category $\hat{\cC}$ whose objects are morphisms in $\cC$ and whose morphisms 
$h_1 \rightarrow h_2$ are commutative diagrams in $\cC$ of the form
\[ \xymatrix{
C_1 \ar[d]^{h_1} \ar[r] & C_2 \ar[d]^{h_2}\\
D_1 \ar[r]              & D_2.
} \]
If no such product exists, $f \e g$ is undefined.  Otherwise, $f \e g$ is well-defined up to isomorphism in $\hat{\cC}$.  

The second operation is only defined for morphisms with the same domain.  Let morphisms $f,g$ have domain $A$.  We
define the \emph{internal product} $f \times_A g \in \hom{(\cC)}$ to be the product, if it exists, of $f$ and $g$ in the category $\cC_A$
whose objects are morphisms with domain $A$, and whose morphisms $h_1 \rightarrow h_2$ are commutative diagrams in $\cC$ of the form
\[ \xymatrix{
A \ar[d]_{h_1} \ar[rd]^{h_2} & \\
C_1 \ar[r] & C_2.
} \]
As above, if no such product exists, $f \times_A g$ is undefined; otherwise, it is well-defined up to isomorphism in $\cC_A$.  

To understand these operations in terms of communication systems, consider the category {\bf FinSet} of finite sets and finite set 
mappings, where each morphism is interpreted as a correspondence between sent and received messages.  
If $f$ and $g$ are set mappings with domain $S$, the internal product $f \times_S g$ is the function $s \mapsto \big (f(s), g(s) \big)$ 
evaluating 
$f$ and $g$ on the same element of $S$.  This corresponds to sending the same message through two different systems simultaneously.  
$f \e g$, on the other hand, is the map $(s_1, s_2) \mapsto \big( f(s_1), g(s_2) \big)$ taking two inputs and evaluating $f$ on the 
first and $g$ on the second.  This corrpesponds to sending two different messages through different systems.  

Somewhat more generally, suppose the product of any two objects exists in our category $\cC$.  Then the external and internal 
products always exist and are, up to isomorphism, the unique morphisms making the following diagrams commute:
\[ \xymatrix{
A \ar[d]^f & A \times B \ar[l] \ar[d]^{f\e g} \ar[r]& B \ar[d]^g\\
X          &     X \times Y  \ar[l] \ar[r]        &     Y
} \]
for external products and
\[ \xymatrix{
              & A \ar[d]|-{f \times_A g} \ar[ld]_f \ar[rd]^g            \\
X      &     X \times Y    \ar[l] \ar[r]         &   Y           
} \]
for internal products.

%%%%%%%%%%%%%%%%%%%%%%%%%%%%%%%%%%%%%%%%%%%%%%%%%%%%%%%%%%%%%%%%%%%%%%%%%%%%%%%%%%%%%%%%%

\section{Axioms for Information Functions}

We now present axioms characterizing the behavior of information measures.  These axioms were chosen because a) they represent 
properties of information we believe to be fundamental, and b) they appear to be a minimal set necessary to establish basic lemmas
(see section 6.)

For a category $\cC$ whose morphisms represent
communication systems, we define an \emph{information function} to be a function that quantifies the amount of information received by 
the destination as a message is sent through a system.  Mathematically, an information function $I$ assigns to each morphism $f$ of
$\cC$ a nonnegative real number $I(f)$, subject to the following axioms. 

First, if two communication systems are mathematically identical, the information flow through each is the same.

\begin{axiom}\emph{Invariance:} If morphisms $f$ and $g$ are isomorphic in $\hat{\cC}$, then $I(f)=I(g)$. \end{axiom}

Second, if two different messages are sent through different systems, so that the received messages are independent of each 
other, the amount of information in both them is the sum of the information in each of them.  

\begin{axiom}\emph{External Additivity:} If $f \e g$ exists, $I(f \e g) = I(f) + I(g)$. \end{axiom}

If the same message is sent through different systems simultaneously, the total information gain cannot exceed the sum of the information
obtained through each system, and may be strictly less due to redundancies in the received messages.  More strongly, we require:

\begin{axiom} \emph{Internal Strong Subadditivity:}  If $f$,$g$, and $h$ have the same domain $A$, then 
\[I(f \times_A g \times_A h) \leq I(f \times_A g) + I(g \times_A h) - I(g), \]
if the above products exist.
\end{axiom}

That is, the information gained from the tripartite system $f \times_A g \times_A h$ cannot exceed the amount gained from 
$f \times_A g$, plus that gained from $g \times_A h$, minus the amount $I(g)$ that is redundant to $f \times_A g$ and $g \times_A h$.

The fourth axiom concerns the case where a single message is sent through two systems in sequence; i.e. the message received from the 
first system is relayed through the second.  No new information can be gained from passing the message through the second system;
this is known as the \emph{data processing inequality}.  Furthermore, if no information is lost in 
sending the message through the second system, then it is possible to reconstruct the first received message from
the second.

\begin{axiom} \emph{Monotonicity:}
Given a diagram
\[ \xymatrix{
A \ar[r]^f & B  \ar[r]^g & C, 
} \]
$I(g \circ f) \leq I(f)$ with equality if and only if there is a morphism $s:C \rightarrow B$ such that $s\circ g \circ f = f$.
\end{axiom}

Finally, the amount of information that can be sent through a system is limited the range of messages which can be received.  
One cannot improve on a system in which the sent messages match the received messages perfectly.

\begin{axiom} \emph{Destination Matching:}  For any morphism $f$ with codomain $B$, $I(f) \leq I(\id_B)$. \end{axiom}

%%%%%%%%%%%%%%%%%%%%%%%%%%%%%%%%%%%%%%%%%%%%%%%%%%%%%%%%%%%%%%%%%%%%%%%%%%%%%%%%%%%%%%%%

\section{Examples of Information Functions}

We now explore mathematical information functions encompassed by our framework, including many classical 
information measures as well as familiar constructs not usually seen as related to information.  We also show how our framwork
might be applied to quantum channel capacities.  With the exception of these quantum capacities, satisfaction of the axioms is
either a well-known result or a simple exercise. 

\subsection{Discrete Communication}

In the simplest mathematical models of communication, there is a finite set of possible messages to send and a 
deterministic mapping from messages sent to messages received.  Such systems can be represented by morphisms in the category
{\bf FinSet} of finite sets and set mappings.

\subsubsection{Hartley Entropy}

The information flow through a discrete system can be quantified in terms of the number of possible received messages.  The greater
variety in the messsages which can be received, the more information is gained upon receiving a particular message.  One simple
way to measure this variety is the following: for a set mapping $f:A \rightarrow B$, we define the \emph{Hartley entropy} of $f$ to be
\[
H_0(f) = \log |f(A)|,
\]
the logarithm (customarily taken base 2) of the cardinality of $f$'s image.  Hartley entropy measures the expected number of bits
needed to encode a received message if all possible received messages are equally likely.

\subsubsection{Shannon Entropy}

The Shannon entropy modifies $H_0$ to incorporate unequal probabilities of received messages.  For a mapping 
$f:A \rightarrow B$ of finite sets, we define the \emph{Shannon entropy} of $f$ to be
\[
H(f) = -\sum_{b \in B} \frac{|f^{-1}(b)|}{|A|} \log \frac{|f^{-1}(b)|}{|A|}.
\]
If all sent messages $a$ are equally likely, $H(f)$ quantifies the uncertainty of the received message as the expected number of bits
needed to encode $f(a)$ using an optimal encoding.  

By a theorem of Aczel, Forte, and Ng \cite{Aczel}, the only information functions on the category {\bf FinSet} are linear 
combinations of the Shannon and Hartley entropies.  

\subsubsection{Noisy Information}

The category {\bf FinSet} has limitations as a setting for modelling discrete communication systems.  There is no way to represent 
communication errors induced by noise, or to record the fact that some messages are more likely to be sent than others.
We can overcome both of these limitations by considering our objects to
have two parts: a visible part representing the messages that can be sent or received,
and a hidden part representing noise that may affect transmission.

Let {\bf NoisyFinSet} be the category whose objects are pairs $(M,A)$ of finite sets, together with a surjective set
map $\pi_A:M \rightarrow A$.  For each $a \in A$, the preimage
$M_a \equiv \pi_A^{-1}(a)\subset M$ represents the environmental noise factors that might be sent along with $a$.  

Morphisms in this category are diagrams
\[ \xymatrix{
M \ar[r]^f \ar[d]^{\pi_A}    &    N \ar[d]^{\pi_B}\\
A                            &    B.
}\]
For each intended message $a \in A$, the actual transmitted data is an element of the preimage 
$M_a$, representing both the original message and the noise.  Different 
elements in this preimage may map to different elements of $B$ under $\pi_B \circ f$, in accordance with the possibility that noise may 
change the received message.  The map $\pi_B$ may be used to model the decoding of an error-correcting code.

For such a morphism we define the \emph{noisy Shannon information} in the following manner: Consider $M$ as a probability space
with normalized counting measure (that is, consider each element of $M$ as equally likely to occur.)  Then $\pi_A$ and 
$\pi_B \circ f$ can be interpreted as random variables with values in $A$ and $B$ respectively.  Define 
$NI(f) = I(\pi_B \circ f; \pi_A)$, the mutual information of these two variables.  This quantity may be computed precisely as
\[
NI(f) = \frac{1}{|M|} \sum_{\substack{a\in A\\b\in B}} |M_a \cap M_b| \log \frac{|M_a \cap M_b|}{|M_b|}- 2 |A| \log |A|,
\]
where $M_b$ is the preimage $f^{-1}(N_b) = f^{-1}(\pi_B^{-1}(b)) \subset M$, i.e. the set of all noise elements $m \in M$ that cause
message $b$ to be received.  This function generalizes the Shannon entropy as defined above.  If an element $m \in M$
is chosen randomly with uniform probability, then $NI(f)$ equals the average amount of Shannon information that 
the received message $\pi_B \circ f(m) \in B$ imparts about the sent message $\pi_A(m)$.   This formalism also
allows for unequal probability in sent messages: the probability of a message $a \in A$ being sent is proportional to the cardinality
of the preimage $M_a$.

\subsubsection{Channel Capacity}

Channel capacity is the maximum amount of information that can be sent through a noisy system,
where the maximum is taken over all probability distributions on the set of sent messages.  To represent channel capacity as
an information function, we again use the categrory {\bf NoisyFinSet}.  Consider a morphism
\[ \xymatrix{
M \ar[r]^f \ar[d]^{\pi_A}    &    N \ar[d]^{\pi_B}\\
A                            &    B.
}\]
For any message $a \in A$ we consider each element of
the preimage $\pi_A^{-1}(a)$ to be equally likely.  Using this rule, given any probability distribution $p_A$ on $A$, we can associate 
a probability distribution on $M$, and, via $\pi_B \circ f$, a joint distribution $p$ on $A \times B$ and a marginal distribution $p_B$
on $B$.  We define the \emph{channel capacity} $C(f)$ to be the maximum over all probability distributions on $A$ of the 
mutual information between the sent message $a \in A$ and the received message $b \in B$:
\[
C(f) = \max_{p_A} \sum_{\substack{a \in A\\b \in B}} p(a,b) \log \frac{p(a,b)}{p_A(a)p_B(b)}.
\]

\subsection{Continuous Communication}

Many real-world communication systems have a continuous range of messages that can be sent or received.   To model these systems 
we use the category \textbf{Prob} whose objects are probability spaces, i.e. measure spaces $(M, \mu)$ with $\mu(M)=1$.   We define
the morphisms in \textbf{Prob} to be measurable functions $f:(M,\mu) \rightarrow (N, \nu)$ that are \emph{backwards measure preserving},
i.e. $\mu(f^{-1}(U)) = \nu(U)$ for each measurable $U \subset N$.  This condition guarantees that the probability measure induced by 
$f$ on $N$ agrees with $\nu$.

Internal products of morphisms are not guaranteed to exist in this category.  The obvious candidate, $x \mapsto (f(x),g(x))$, 
does not in general satisfy the backwards measure preserving condition.  External products, however, do exist: for 
$f:(M_1, \mu_1) \rightarrow (M_2, \mu_2)$ and $g:(N_1, \nu_1) \rightarrow (N_2, \nu_2)$, the external product $f \e g$ is the 
natural map between the product spaces $M_1\times N_1$ and $M_2 \times N_2$, sending $(x,y)$ to $(f(x), g(y))$.  It is easily
verified that this map is backwards measure preserving.

\subsubsection{Noisy Information}

To represent imperfect communication, we add ``noise'' to this category in the same manner as for \textbf{FinSet}, by forming a 
category \textbf{NoisyProb} whose objects are pairs $\big( (M,\mu), (A, \alpha) \big)$ with morphisms $\pi_A:M \rightarrow A$.  
Morphisms in \textbf{NoisyProb} are again given by commutative diagrams
\[ \xymatrix{
(M,\mu) \ar[r]^f \ar[d]^{\pi_A}    &    (N,\nu) \ar[d]^{\pi_B}\\
(A, \alpha)                        &    (B, \beta).
}\]

Such a morphism $f$ induces a surjective measurable map $\tilde{f}: M \rightarrow A \times B$ sending $m \in M$ to $(\pi_A(m),
\pi_B \circ f(m))$.    
$\tilde{f}$ and the measure $\mu$ on $M$ induce a probability measure $\rho$ on $A \times B$.  If $\rho$ is absolutely
continuous with respect to $\mu \times \nu$, we define the \emph{(continuous) noisy Shannon information} of $f$ to be
\[
NI_{cont}(f) = \int_{A \times B} \ln \frac{d \rho}{d (\mu \times \nu)} \; d\rho,
\]
where $\frac{d \rho}{d (\mu \times \nu)}$ is the Radon-Nikodym derivative of $\rho$ with respect to $\mu \times \nu$.
If $\rho$ is not absolutely continuous with respect to $\mu \times \nu$, we leave $I(f)$ undefined.  
Our definition agrees with the usual formula for continuous mutual information
\[
I = \int_{A \times B} p(x,y) \ln \frac{p(x,y)}{p_A(x)p_B(y)} \; d(\mu \times \nu)
\]
under the substitutions $p(x,y) = \frac{d \rho}{d (\mu \times \nu)}$ 
and $p_A(x) \equiv p_B(y) \equiv 1$ with respect to the probability measures $\alpha$ and $\beta$, respectively.

(It may come as a surprise that $NI_{cont}$ for a \emph{noiseless} system, wherein $\pi_A$ is an isomorphism, is usually 
undefined if $A$ contains uncountably many points.  This is because a communication system that
can faithfully transmit any of uncountably many inputs can send any finite amount of information in a single message.)

\subsubsection{Channel Capacity}

Channel capacity is obtained by maximizing the amount of transmitted information over the set of all probability distributions on the
space of sent messages.  We represent such probability distributions by measurable functions $p: A \rightarrow \Rp$ with
$\int_A p \; d\alpha = 1$.  For a particular message $a \in A$, we view all environmental noise factors in the preimage
$\pi_A^{-1}(a)$ as equally likley; we therefore define a new probability measure $\mu_p$ on $M$ by
\[
\mu_p (U) = \int_U p\circ \pi_A \; d\mu
\]
for each measurable $U \subset M$.  

Given a morphism $f$ as above, the map $\tilde{f}:M \rightarrow A \times B$ induces a probability measure $\rho_p$ on 
$A\times B$.  The \emph{continuous channel capacity} of $f$ is defined as
\[
C_{cont}(f) = \sup_p \int_{A \times B} \ln \frac{d \rho_p}{d (\mu \times \nu)} \; d\rho_p,
\]
if such a supremum exists; otherwise $C(f)$ is undefined.

\subsection{Quantum Communication}

For systems that communicate using processes governed by quantum mechanics, the transmission of information is described by quantum
information theory.  
Quantum information theory is a young field with many open problems; in particular, we cannot say at this time whether 
the various capacities of a quantum channel satisfy our axioms for information functions.  However, our category-theoretic
framework suggests new problems and casts existing ones in a new light.

Consider a finite-dimensional quantum system represented by Hilbert space $H$.  A ``message'' in this system is represented by 
a density matrix: a positive-semidefinite Hermitian operator of trace one acting on $H$.  Let $D(H)$ denote
the space of all density matrices over the Hilbert space $H$.  Given two such spaces $D(H_1)$ and $D(H_2)$, the natural
choice for a morphism between them is a completely positive trace-preserving linear map.  We call such maps \emph{CP maps}
or \emph{quantum channels}.  We denote by \textbf{Quant} the category whose objects are spaces $D(H)$ of density 
matrices over finite-dimensional Hilbert spaces, and whose morphisms are CP maps.

External products in this category are given by tensor product of CP maps.  Internal products do not exist in all 
cases;  for instance, the internal product
of the identity map $\rho \mapsto \rho$ with itself would have to be $\rho \mapsto \rho \otimes \rho$, but this is nonlinear and
violates the no-cloning theorem.  However, there are important cases where internal products do exist; for example, given a
bipartite system $AB$, the internal product of the partial traces $\Tr_A:D(A \otimes B) \rightarrow D(A)$ and 
$\Tr_B:D(A \otimes B) \rightarrow D(A)$ is the identity map $D(A \otimes B) \rightarrow D(A \otimes B)$.  

There are several different notions of capacity for a quantum channel, depending on the intended application.  The
\emph{quantum capacity} measures the capacity of the channel to send quantum states intact.  The \emph{classical capacity} 
measures the capacity to send classical information.  The \emph{entanglement-assisted classical capacity} is the capacity
to send classical information if the sender and receiver are allowed to share an arbitrary number of entangled quantum states
prior to transmission.  For mathematical definitions of these quantities we refer our readers to the literature \cite{Nielsen, 
Keyl}.

All the above capacities satisfy invariance, monotonicity, and destination matching.  External additivity for these capacities
is a famous open problem; it amounts to the question of whether quantum entangled messages can be communicated more efficiently
than unentangled messages.  For partial results on this problem, see \cite{King, Shor}.
Internal strong subadditivity of these capacities has not, to our knowledge, been studied.  This
question is likely related to the strong subadditivity of Von Neumann entropy, proven by Lieb and Ruskai \cite{Lieb}.

\subsection{Vector Communication}

There are strong reasons to suspect a relationship between information and dimension.  For example, the dimension of a 
vector space equals the number of coordinates
required to specify a point in the space, much like the Shannon entropy equals the average number of bits needed to specify the value 
of a random variable.  Deep connections have also been found between notions of dimension and complexity \cite{Lutz, Hitchcock}.  It 
is therefore not surprising that vector space dimension is also an information function.

Suppose we have a communication system in which the set of possible messages is represented by a finite-dimensional vector space, and
sent messages are mapped linearly onto received messages.  
We represent such systems using the category $K$-$\mathbf{FinVect}$ of finite-dimensional vector spaces over a field $K$, with 
morphisms given by linear maps.  For a morphism $f:V \rightarrow W$, we define
\[
d(f) = \dim f(V).
\]
$d(f)$ can be interpreted as the amount by which knowledge of $f(v) \in W$ reduces the dimensional uncertainty of $v \in V$.

\subsection{Information Functions on Dual Categories}
\label{Dual Categories}

Finally, we present a pair of information functions whose interpretation in terms of communication systems is not clear.  The
categories for these are the duals of the familiar categories \textbf{FinSet} and $K$-$\mathbf{FinVect}$.   In a dual category, 
the internal and external product operations correspond to similarly defined internal and
external \emph{coproduct} operations in the original category.

\subsubsection{Cardinality of Image}

Consider the category $\mathbf{FinSet}^*$ whose objects are finite sets and whose morphisms $f^*:A \rightarrow B$ correspond to
set mappings $f:B \rightarrow A$.  The external product of $f^*:A \rightarrow B$ and $g^*:C \rightarrow D$ is the morphism
$f^* \e g^*:A \sqcup C \rightarrow B \sqcup D$ ($\sqcup$ denotes disjoint union) corresponding to the set mapping
\begin{align*}
f \hat{+} g & : B \sqcup D  \rightarrow A \sqcup C\\
f \hat{+} g & : x \mapsto \begin{cases}
f(x) & \text{if } x \in B\\
g(x) & \text{if } x \in D.
\end{cases}
\end{align*}
The internal product of $f^*:A \rightarrow B$ and $g^*:A \rightarrow C$ is the morphism
$f^* \e g^*:A \rightarrow B \sqcup C$ corresponding to the set mapping
\begin{align*}
f +_A g & : B \sqcup C  \rightarrow A\\
f +_A g & : x \mapsto \begin{cases}
f(x) & \text{ if } x \in B\\
g(x) & \text{ if } x \in C.
\end{cases}
\end{align*}

The function mapping $f^*$ to $|f(B)|$, the cardinality of the image of the corresponding set map $f$, is an information 
function in this category.  It is unclear, however, in what sense the dual of a set mapping represents a communication
system.

\subsubsection{Dimension of Image}

Similarly, we can explore the category $K$-$\mathbf{FinVect}^*$ whose objects are finite 
dimensional vector spaces and whose morphisms $f^*:V \rightarrow W$ correspond to
linear maps $f:W \rightarrow V$.  The external product of $f^*:V \rightarrow W$ and $g^*:T \rightarrow U$ is the morphism
$f^* \e g^*:V \oplus T \rightarrow W \oplus U$ corresponding to the linear map
\begin{align*}
f \hat{+} g & : W \oplus U  \rightarrow V \oplus T\\
f \hat{+} g & : w \oplus u \mapsto f(w) \oplus g(u) \quad \text{for $w \in W, u \in U$}.
\end{align*}
The internal product of of $f^*:V \rightarrow W$ and $g^*:V \rightarrow U$ is the morphism
$f^* \e g^*:V \rightarrow W \oplus U$ corresponding to the linear map
\begin{align*}
f +_V g & : W \oplus U  \rightarrow V\\
f +_V g & : w \oplus u \mapsto f(w) + g(u) \quad \text{for $w \in W, u \in U$}.
\end{align*}

%%%%Fix the following paragraph!!!!

For a morphism $f^*:V \rightarrow W$, the function $I(f^*) = \dim(f(W))$ is an information function in this category.
Since dimension of image is also a information function in the original cageory $K$-$\mathbf{FinVect}$, this function
might be called a \emph{bi-information function}, in that it is an information function on a category and its dual.  The 
existence of a bi-information function for $K$-$\mathbf{FinVect}$ is doubtless related to the existence of biproducts in 
this category.

%%%%%%%%%%%%%%%%%%%%%%%%%%%%%%%%%%%%%%%%%%%%%%%%%%%%%%%%%%%%%%%%%%%%%%%%%%%%%%%%%%%%%%%%%%%%%%

\section{Basic Results}

We now establish some basic facts about information functions.  We start by showing that products and compositions are 
well-defined on isomorphism classes in $\hat{\cC}$ and $\cC_A$.

\begin{proposition} 
\label{isomorphisms}
Suppose $a \cong b$ and $c \cong d$ in $\hat{\cC}$, and $e \cong f$ and $g \cong h$ in $\cC_A$.  Then, assuming
the following compositions and products exist,
\begin{enumerate}
\item $a \circ c \cong b \circ d$ in $\hat{\cC}$.
\item $a \e c \cong b \e d$ in $\hat{\cC}$.
\item $a \circ e \cong b \circ f$ in $\cC_A$.
\item $e \times_A g \cong f \times_A h$ in $\cC_A$.
\end{enumerate}
\end{proposition}

%%Is it worth stating commutativity and associativity?

\begin{proof} Exercise. \end{proof}

Now let $I$ be an arbitrary information function.  The following proposition says that no communication system can do
better than one which reproduces the inputs exactly.
\begin{proposition}
\emph{Source matching:} For any morphism $f$ with domain $A$, $I(f) \leq I(\id_A)$.
\end{proposition}

\begin{proof} $I(f) = I(f \circ \id_A) \leq I(\id_A)$ by monotonicity. \end{proof}

\begin{proposition}
\emph{Monotonicity with respect to internal products:}  Let $f$ and $g$ be morphisms with domain $A$.  Assuming the relevant products 
exist,
\begin{enumerate}
\item $I(f \times_A g) \geq I(f)$
\item $I(f \times_A f) = I(f).$
\end{enumerate}
\end{proposition}

In other words, you will never lose any information by sending the same message simultaneously through two different systems 
versus sending it through just one of them, but if the two systems are identical, you won't gain any information either.

\begin{proof}
\begin{enumerate}
\item By the definition of products there is a morphism $\pi_1:f \times_A g \rightarrow f$ in $\cC_A$, which corresponds to
a diagram
\[ \xymatrix{
A \ar[d]_{f \times_A g} \ar[rd]^{f} & \\
C \ar[r] & B.
} \]
The result now follows from monotonicity.
\item We must now show $I(f \times_A f) \leq I(f)$.  Again by the definition of products, there is a unique morphism 
$\delta:f \rightarrow f \times f$ in $\cC_A$ making the following diagram commute:
\[ \xymatrix{
       & f \ar[d]^{\delta} \ar[ld]_{\id} \ar[rd]^{\id}            \\
f      &     f \times_A f    \ar[l] \ar[r]                  &   f.           
} \]
The morphism $\delta$ represents a diagram
\[ \xymatrix{
A \ar[d]_{f} \ar[rd]^{f \times_A f} & \\
B \ar[r] & C,
} \]
from which the result follows by monotonicity. \qedhere
\end{enumerate}
\end{proof}

Further results depend on the existence of terminal objects in our category $\cC$.  A object $T$ is terminal if there
is exactly one morphism into $T$ from any object $A$.  Any singleton set is a terminal object in \textbf{FinSet}.  Intuitively,
a destination represented by a terminal object is unable to discriminate between messages, and thus cannot receive information.
We prove this and other results in the following propositions:

\begin{proposition} 
\label{terminal1}
Let $T$ be a terminal element of $\cC$.  For a fixed object $A$, let $t$ be the unique morphism 
$t: A \rightarrow T$.  Then, assuming the appropriate products exist,
\begin{enumerate}
\item  $\id_T$ and $t$ are terminal elements in $\hat{\cC}$ and $\cC_A$, respectively.
\item  For any morphism $g \in \hom\cC$, $g \e \id_T \cong g$ in $\hat{\cC}$.  If $g$ has domain $A$, then 
$g \times_A t \cong g$ in $\cC_A$.
\item  Let $B$ be an object for which $A \times B$ exists.  Then $\pi_B \cong t \e \id_B$ in $\cC_{A\times B}$.
\end{enumerate}
\end{proposition}

\begin{proof}
\begin{enumerate}
\item Exercise.
\item Combine part (a) with the fact that $A \times T \cong A$ for any object $A$, in any category with terminal object $T$.
\item Let $s$ and $t$ be the unique morphisms $t:A \rightarrow T$, $s:B \rightarrow T$.  Let $e$ be the unique morphism 
$e:B \rightarrow T \times B$ making the following diagram commute:
\[ \xymatrix{
&  B \ar[ld]_s \ar[d]^e \ar[rd]^{\id_B}\\
T & T\times B \ar[l] \ar[r] & B
} \]
Then the diagrams
\[ \xymatrix{
& A \times B \ar[ld]_{t \circ \pi_A} \ar[d]|-{e \circ \pi_B} \ar[rd]^{\pi_B}\\
T & T\times B \ar[l] \ar[r] & B
} \]
and
\[ \xymatrix{
& A \times B \ar[ld]_{t \circ \pi_A} \ar[d]|-{t \e \id_B} \ar[rd]^{\pi_B}\\
T & T\times B \ar[l] \ar[r] & B
} \]
also commute.  By the definition of products we must have $e \circ \pi_B = t \e \id_B$, and the result now follows from the
fact that $e$ is an isomorphism. \qedhere
\end{enumerate}
\end{proof}

\begin{corollary}
Let $f: A \rightarrow C$ be any morphism, and let $e$ and $t$ be as above.  Then
\begin{enumerate}
\item $I(f \e \id_T) = I(f \times_A t) = I(f)$ for any morphism $f$, whenever these products are defined.
\item $I(\id_T) = I(t) = 0$.
\item \emph{Subadditivity:} $I(f \times_A g) \leq I(f) + I(g)$ for any two morphisms $f,g$ with domain $A$.
\item If the product $A\times B$ exists for an object $B$, then $I(f \circ \pi_A) = I(f)$.
\end{enumerate}
\end{corollary}

Part (d) says that the information flow through a communication system $f$ is not affected by the presence of ``irrelevant''
information in $B$.

\begin{proof}
Parts (a) follows direcly from the proposition.  For part (b), $I(\id_T)=0$ follows from part (a) and additivity, then $I(t)=0$ 
by destination matching.  Part (c) is proven by subsituting $t$ into strong subadditivity,
\[I(f \times_A t \times_A g) + I(t) \leq I(f \times_A t) + I(t \times_A g), \]
and invoking previous results.  For part (d), $I(\pi_A) = I(\id_A)$ by Proposition 
\ref{terminal1}(c), and the result follows from Proposition \ref{isomorphisms}(a).
\end{proof}

\section{Conclusion}

Mathematical abstractions like groups and topological spaces have the power to illuminate connections between different
objects of study and inspire the discovery of new objects.  We hope this present abstraction will lead to
the discovery of new ways to measure information, and deepen understanding of information's general properties.
In the future, we aim to prove deeper results and investigate information functions relevant to specific applied situations.  

\bibliography{uinfo.bib}{}
\bibliographystyle{plain}

\end{document}